\documentclass[12pt,oneside]{article}
\ifdefined\HCode   

\fi
\usepackage{geometry}
\usepackage{graphicx}
\usepackage{amsmath}
\usepackage{amsfonts}
\usepackage{amssymb}
\usepackage{tikz-cd}
\usepackage{comment}
\usepackage{fullpage}
\usepackage{epstopdf}
\usepackage{amsthm}
\usepackage{mathrsfs} 
\usepackage{tikz}
\usepackage{dsfont}
\usetikzlibrary{matrix}
\usepackage{mathtools}
\usepackage[hidelinks]{hyperref}
\usepackage{cleveref}
\usepackage{xcolor}
\usepackage{enumitem}
\usepackage{tensor}
\usepackage[utf8]{inputenc}
\usepackage{csquotes}
\usepackage[english]{babel}
\usepackage[
backend=bibtex,
style=alphabetic,
sorting=nyt
]{biblatex}
\usepackage{mymacros}

\hypersetup{
	colorlinks,
	linkcolor={red!50!black},
	citecolor={blue!50!black},
	urlcolor={blue!80!black}
}

\tikzset{
	rot90/.style={anchor=south, rotate=90, inner sep=.5mm}
}
\tikzset{
	rot45/.style={anchor=south, rotate=-45, inner sep=.5mm}
}

\newtheorem{theorem}{Theorem}[section]

\newtheorem{lemma}[theorem]{Lemma}
\newtheorem{definition}[theorem]{Definition}
\newtheorem{proposition}[theorem]{Proposition}
\newtheorem{corollary}[theorem]{Corollary}

\theoremstyle{definition}
\newtheorem{construction}{Construction}[theorem]
\newtheorem{question}{Question}[theorem]
\newtheorem{remark}{Remark}[theorem]

\newtheorem{innercustomthm}{Theorem}
\newenvironment{customthm}[1]
{\renewcommand\theinnercustomthm{#1}\innercustomthm}
{\endinnercustomthm}

\theoremstyle{definition}

\newtheorem*{dfn*}{Definition}
\newtheorem*{axm*}{Axiom}
\newtheorem*{ntn*}{Notation}
\newtheorem*{exm*}{Example}
\newtheorem*{exr*}{Exercise}
\newtheorem*{int*}{Intuition}
\newtheorem*{qst*}{Question}
\newtheorem*{rmk*}{Remark}

\theoremstyle{plain}

\newtheorem*{thm*}{Theorem}
\newtheorem*{prop*}{Proposition}
\newtheorem*{cor*}{Corollary}
\newtheorem*{lem*}{Lemma}
\newtheorem*{cnj*}{Conjecture}

\title{Eilenberg Mac Lane spectra as p-cyclonic Thom spectra}
\author{Ishan Levy\thanks{The author is supported by the NSF Graduate Research Fellowship under Grant No. 1745302.}}
\usepackage{setspace}

\makeatletter
\renewcommand\tableofcontents{%
	\@starttoc{toc}%
}
\makeatother

\begin{document}
	\date{}
	\maketitle
	\begin{abstract}
		Hopkins and Mahowald gave a simple description of the mod $p$ Eilenberg Mac Lane spectrum ${\FF}_p$ as the free $\EE_2$-algebra with an equivalence of $p$ and $0$. We show for each faithful $2$-dimensional representation $\lambda$ of a $p$-group $G$ that the $G$-equivariant Eilenberg Mac Lane spectrum $\ul{\FF}_p$ is the free $\EE_{\lambda}$-algebra with an equivalence of $p$ and $0$. This unifies and simplifies recent work of Behrens, Hahn, and Wilson, and extends it to include the dihedral $2$-subgroups of $\Or(2)$. The main new idea is that $\ul{\FF}_p$ has a simple description as a $p$-cyclonic module over $\THH(\FF_p)$. We show our result is the best possible one in that it gives all groups $G$ and representations $V$ such that $\ul{\FF}_p$ is the free $\EE_V$-algebra with an equivalence of $p$ and $0$.
	\end{abstract}
	\begin{spacing}{0.1}
		\tableofcontents
	\end{spacing}
	\section{Introduction}
	
	Consider the nontrivial virtual bundle on the circle, $S^1 \to \BO$. We can extend this to a double loop map $\mu_2: \Omega^2S^3 \to \BO$. Mahowald proved the following result \cite{MAHOWALD1977249}:
	
	\begin{theorem}[Mahowald]\label{mahowald}
		The Thom spectrum $(\Omega^2S^3)^{\mu_2}$ is equivalent as an $\EE_2$-ring to the Eilenberg Mac Lane spectrum $\FF_2$.
	\end{theorem}
	
	This result and its generalizations have since found great use, for example giving a conceptual explanation of Bokstedt's computation of $\THH(\FF_p)$ \cite{blumberg2008thh,blumberg2010topological}.
	
	Hopkins observed that Mahowald's theorem generalizes to odd primes. We can consider the $p$-local stable spherical fibration on the circle $S^1 \to \BGL_1(\SP_{(p)})$ representing the class $1-p$ in $\pi_1$. Extending this to a double loop map $\mu_p:\Omega^2S^3 \to \BGL_1(\SP_{(p)})$, he showed:
	
	\begin{theorem}[Hopkins {\cite[Theorem 4.18]{Mathew_2015}}]\label{hopkins}
		The Thom spectrum $(\Omega^2S^3)^{\mu_p}$ is equivalent as an $\EE_2$-ring to the Eilenberg Mac Lane spectrum $\FF_p$.
	\end{theorem}
	
	\begin{remark}\label{freeinterpretation}
		The Thom spectrum in Theorems \ref{mahowald},\ref{hopkins} can be interpreted as the free $\EE_2$-ring with $p \simeq 0$ \cite[Theorem 5.1]{simpleunivpropofthomspectra}. This means an $\EE_{2}$-algebra map $(\Omega^2S^3)^{\mu_p}  \to R$ is the same as a nulhomotopy of $p$ in $R$.
	\end{remark}
	
	Behrens--Wilson proved an equivariant generalization of Mahowald's result \cite{behrens2018c2equivariant}, which we now explain.
	
	\begin{construction}\label{loopmap}
		Let $G$ be a group and $V: G \to \Or(n)$ be any representation. We can consider the $V$-loop map $\mu_{V,p}:\Omega^{V}S^{1+ V} \to \BGL_1(\SP_{(p)})$ extending the map $1-p \in \pi^G_1(\BGL_1(\SP_{(p)}))$ where $\SP$ is the $G$-equivariant sphere spectrum. Then $X_{V,p}$ is defined to be the Thom spectrum of $\mu_{V,p}$, which is an $\EE_V$-algebra.
	\end{construction}
	
	The Thom spectra $X_{V,p}$ are equivariant generalizations of the Thom spectra considered by Hopkins and Mahowald. For example, as in \Cref{freeinterpretation}, $X_{V,p}$ is the free $\EE_{V}$-algebra with $p\simeq 0$.
	
	\begin{theorem}[Behrens--Wilson]\label{bw}
		Let $\rho$ be the regular representation of the group $C_2$. Then there is an equivalence of $\EE_{\rho}$-algebras $X_{\rho,2} \cong \ul{\FF}_2$.
	\end{theorem}
	
	Their proof relies on studying $\EE_\rho$ power operations, similarly to how Mahowald and Hopkins' results rely on  $\EE_2$ power operations.
	
	Hahn--Wilson generalized the result to all cyclic $p$-subgroups of $\Or(2)$ \cite{hahn2019eilenbergmaclane}.
	\begin{theorem}[Hahn--Wilson]\label{hw}
		Let $\lambda:G \to \Or(2)$ be a faithful representation of a cyclic $p$-group. Then there is an equivalence of $\EE_{\lambda}$-algebras $X_{\lambda,p}\cong \ul{\FF}_p$.
	\end{theorem}
	
	Notably, their proof did not involve studying $\EE_\lambda$ power operations, but rather used the structure on $X_{\lambda,p}$ as a module over its norm $N_*^G(X_{\lambda,p})$, the nonequivariant Hopkins--Mahowald result, and an exotic Thom ring structure on the Thom spectrum as inputs.
	
	We extend the result of Hahn--Wilson to the dihedral $2$-subgroups of $\Or(2)$ and give a simpler proof in the case of $C_{p^n}\subset \SOr(2)$, removing the need for exotic Thom ring structures.
	
	\begin{customthm}{A} \label{mainthm}
		Let $\lambda:G \to \Or(2)$ be a faithful representation of a $p$-group. Then there is an equivalence of $\EE_{\lambda}$-algebras $X_{\lambda,p}\cong \ul{\FF}_p$.
	\end{customthm}
	There is a concise way to capture all of the structure present in \Cref{mainthm}. Recall \cite[Definition 2.10]{sulyma2020slice} that a \textit{cyclonic} spectrum is one that has a Borel $S^1$-action genuine with respect to the finite subgroups of $S^1$\footnote{This terminology was introduced in \cite{barwick2016cyclonic}.}. If it is only genuine with respect to the finite $p$-subgroups of $S^1$, it is said to be \textit{p-cyclonic}.
	 
	 In this language, for $p>2$, \Cref{mainthm} says that $X_{\lambda,p}\cong \ul{\FF}_p$ as $p$-cyclonic $\EE_{\lambda}$-algebras\footnote{The equivalence is a map of genuine $\Or(2)$ spectra (though it is only an equivalence on the $p$-subgroups).}.
	 
	 At the prime $2$, \Cref{mainthm} gives more than a $2$-cyclonic equivalence since it is a genuine equivalence on the finite dihedral  $2$-groups $D_{2^{n}}$ as well. Accordingly, we can say that a \textit{real 2-cyclonic} spectrum is a Borel $\Or(2)$-spectrum that is genuine on the finite $2$-subgroups of $\Or(2)$. Then \Cref{mainthm} says for $p=2$ that $X_{\lambda,2}\cong \ul{\FF}_2$ as real $2$-cyclonic $\EE_{\lambda}$-algebras.
	 \vspace{7pt}
	
	\Cref{mainthm} answers a question of Hahn and Wilson \cite[Question 9.4]{hahn2019eilenbergmaclane}, asking if $\ul{\FF}_p$ can be a Thom spectrum for any group $G$ that is not a cyclic $p$-group. It is also the best possible result in that it gives all groups $G$ and representations $V$ of $G$ such that the free $\EE_{V}$-algebra with $p \simeq 0$ is $\ul{\FF}_p$.
	
	\begin{customthm}{B} \label{nogothm}
		Let $G$ be a finite group and $V:G\to \Or(n)$ a representation. $X_{V,p} \cong \ul{\FF}_p$ if and only if $n = 2$, $G$ is a p-group, and $V$ is faithful.
	\end{customthm}
	
	The main new idea in the proof of \Cref{mainthm} is to use a module structure that exists on both $\ul{\FF}_p$ and $X_{\lambda,p}$ that is constructed via equivariant factorization homology.
	
	\begin{definition}
		Let $G$ be a finite group and $V:G\to \Or(n)$ a representation. For an $\EE_{V}$-algebra $X$, we define $N^{V}X$ to be the genuine equivariant factorization homology $\int_{\RR^{V}-\RR^{V^G}}X$ as in 
		\cite{horev2019genuine}.
	\end{definition}
	
	Some important properties of $N^{V}$ are:
	\begin{enumerate}[label = (N\arabic*)]
		\item Any $\EE_{V}$-algebra $R$ is naturally an $\EE_{V^G}$-algebra over the $\EE_{V^G+1}$-ring $N^{V}R$.
		\item For any map $f$ of $\EE_V$-algebras such that $f^{\Phi H}$ is an equivalence for any proper subgroup $H\subset G$, $N^{V}f$ is an equivalence.
	\end{enumerate}
	
	In the proof of \Cref{mainthm}, we construct a map of $\Or(2)$-$\EE_{\lambda}$-algebras $f:X_{\lambda,p} \to \ul{\FF}_p$, and prove by induction on the finite $p$-subgroups $G \subset \Or(2)$ that $f^{\Phi G}$ is an equivalence. By (N2), in the induction step we know $(N^{\lambda}f)^{\Phi G}$ is an equivalence, so that by (N1), $f^{\phi G}$ is a $N^{\lambda}\ul{\FF}_p$-module map. The following result is key in showing $f^{\Phi G}$ is an equivalence.
	
	\begin{customthm}{C}\label{freemodrank2both}
		Let $\lambda:G \to \Or(2)$ be a faithful irreducible representation of a $p$-group. Then $X_{\lambda,p}^{\Phi G}$ and $\ul{\FF}_p^{\Phi G}$ are rank $2$ free modules over $(N^{\lambda}X_{\lambda,p})^{\Phi G}$ and $(N^{\lambda}\ul{\FF}_p)^{\Phi G}$ respectively, generated in degrees $0$ and $1$.
	\end{customthm}
	
	Because $f^{\Phi G}$ is a module map, \Cref{freemodrank2both} reduces showing $f^{\Phi G}$ is an equivalence to just showing it is $1$-connected.
	
	Our proof of \Cref{freemodrank2both} is a calculation, which for the cyclic groups $G=C_{p^n}$ is made easier by the fact that we identify $N^{\lambda}\ul{\FF}_p$ as a $G$-spectrum with the topological Hochschild homology $\THH(\FF_p)$. In fact, for the dihedral $2$-subgroups $\lambda:D_{2^n} \to \Or(2)$, $N^{\lambda}\ul{\FF}_2$ agrees with the real topological Hochschild homology $\THR(\FF_2)$.
	
	\begin{remark}
		\Cref{freemodrank2both} can be also be concisely stated for all $G$ at once using cyclonic spectra. At odd primes, it says that $X_{\lambda,p}^{\Phi C_p}$ and $\ul{\FF}_p^{\Phi C_p}$ are free rank $2$ $p$-cyclonic $\THH(\FF_p)$-modules in degrees $0$ and $1$, and for $p=2$, the statement is the same, but with real $2$-cyclonic instead of $p$-cyclonic.
	\end{remark}
	
	In \cref{other}, we also identify $X_{\lambda,p}$ when $\lambda$ is a faithful representation of a cyclic group $C_q$ for $q$ a prime not equal to $p$. The underlying $G$-spectrum is algebraic: it is a sum of suspensions of Mackey functors. More generally, we can ask for $X_{\lambda,p}$ to be algebraic in a weaker sense: the underlying $G$-spectrum, viewed as a spectral Mackey functor as in \cite{guillou2017models} or \cite{barwick2014spectral}, can be lifted to a Mackey functor with values in $\FF_p$-modules.
	
	\begin{question}\label{conjecture}
		Under what conditions does the $G$-spectrum $X_{V,p}$ lift to a Mackey functor with values in $\FF_p$-modules?
	\end{question}
	A necessary condition is that $\dim(V)\geq 2$, and in \cref{other}, we explain why it is plausible this is sufficient.
	
	In \cref{deloop}, we show that when $p=2$, the Thom spectra in \Cref{mainthm} admit extra structure.
	
	\begin{customthm}{D}\label{deloopthm}
		Let $\lambda$ be the standard $2$-dimensional representation of $\Or(2)$. The map $\mu_{\lambda,2}:\Omega^{\lambda}S^{1+ \lambda} \to \BGL_1(\SP_{(2)})$ admits a refinement to a $\lambda+\det\lambda$-loop map.
	\end{customthm}
	
	\begin{corollary}\label{extrastr}
		When $p=2$, the equivalence in \Cref{mainthm} is one of $\EE_{\lambda+\det\lambda}$-algebras.
	\end{corollary}
	
	It is known that \Cref{mahowald} refines to an $\EE_{3}$-equivalence, coming from the equivalence $\Omega^2S^3 = \Omega^3 \BSp(1)$. In \cite{hahn2019eilenbergmaclane}, it was shown that the equivalence in \Cref{hw} is an $\EE_{\lambda+1}$ equivalence for $G$ a cyclic $2$-subgroup of $\Or(2)$. \Cref{extrastr} extends these results to include the dihedral subgroups $D_{2^n}$ of $\Or(2)$.
	
	\subsection*{Acknowledgements}
	
	I am very grateful to Jeremy Hahn for suggesting this project to me and for many discussions about it. I am also grateful to Robert Burklund, Mike Hopkins, Asaf Horev, Dylan Wilson, and Allen Yuan for discussions related to this project. I also thank Robert Burklund and Jeremy Hahn for helpful comments on earlier drafts. Finally, I thank the referee for their useful feedback.

	\section{Module structures for $\EE_V$-algebras}\label{modulestrs}
	
	An important tool in the proof of \Cref{mainthm} is the use of various module and algebra stuctures associated to $\EE_V$-algebras. In this section, we explain how these structures can be obtained via genuine equivariant factorization homology as in \cite{horev2019genuine}.
	
	Equivariant factorization homology is a construction taking as input a $V$-framed $G$-manifold $M$ and an $\EE_{V}$-algebra $X$, and giving an object $\int_{M}X$ as output. For definitions and details, see \cite{horev2019genuine}, but we recall some properties here:
	
	\begin{enumerate}[label = (F\arabic*)]
		\item $\int_{\RR^{V}}X = X$.
		\item It is a $G$-symmetric monoidal functor, meaning it is natural with respect to restriction to subgroups and takes $G$-set indexed disjoint unions of manifolds to $G$-set indexed tensor products (i.e the Hill-Hopkins-Ravenel norm in $G$-spectra). 
		\item If a $G$-framed manifold $M$ is equipped with a collar, that is an identification $M \simeq M'\otimes (-1,1)$ as a framed $G$-manifold, this naturally gives $\int_{M'\otimes (-1,1)}X$ an $\EE_1$-algebra structure \cite[Construction 5.2.1]{horev2019genuine}.
		\item ($G$-$\otimes$-excision) Given a $V$-framed manifold $M$ with a decomposition $M = M_{-1} \cup_{N\times (-1,1)}M_1$ into open submanifolds $M_{-1},M_1$ glued along a collar $N\times (-1,1)$,  $\int_{M_{-1}}X$ and $\int_{M_{1}}X$ are naturally left and right modules respectively over $\int_{N\times(-1,1)}X$, and $\int_M X = \int_{M_{-1}}X\otimes_{\int_{N\times(-1,1)}X}\int_{M_{1}}X$.
	\end{enumerate}
	
	From these properties, we extract general results about $\EE_{V}$-algebras.
	\begin{proposition}\label{mod}
		Let $U$ be an open subset of $S(V)$, the unit sphere of $\RR^V$, and $X$ an $\EE_{V}$-algebra. Let $S_U$ be the preimage of $U$ under the projection $\RR^V-0 \to S(V)$. Then $X$ is a module over the $\EE_1$-ring $\int_{S_U}X$, and there is a canonical map of $\int_{S_U}X$ modules $\int_{S_U}X \to X$ that is compatible with unit maps $\SP \to X$ and $\SP \to \int_{S_U}X$. This module structure is natural in $U$.
		\begin{proof}
			The $V$-framed $G$-manifold $\RR^{V}$ decomposes as $\RR^{V}\cup_{S_U} S_U$, and $S_U$ admits a collar $S_U\cong U\times (-1.1)$. Thus the first claim follows from property $(F4)$ of equivariant factorization homology above. The second claim follows by observing that there is a tautological decomposition $S_U = S_U\cup_{S_U} S_U$ compatible with the decomposition of $\RR^{V}-0$, so by naturality of property $(F4)$, we obtain our map of modules. To see that this map is unital, we simply observe that the unit maps are the maps induced on factorization homology from the inclusion of the empty manifold $\phi$. Naturality in $U$ follows from naturality in $(F3),(F4)$.
		\end{proof}
	\end{proposition}
	
	When we set $U = S(V)-S(V^G)$ in \Cref{mod} $S_U = \RR^V-\RR^{V^G}$, so we obtain:
	\begin{corollary}\label{snormmod}
		Given an $\EE_V$-algebra $X$, $X$ is naturally a module over $N^VX$, and receives a unital $N^VX$-module map from $N^VX$.
	\end{corollary}
	
	\begin{remark}
		\Cref{snormmod} says $X$ is an $\EE_0$-algebra in $N^VX$-modules. In fact, it is true that $N^VX$ is an $\EE_{V^G+1}$-algebra and $X$ is canonically an $\EE_{V^G}$-algebra in $N^VX$-modules, and this refines the $\EE_{V^G}$-algebra in $G$-spectra.
	\end{remark}
	
	We also need the following module structure, where $N_*^PX$ is the indexed smash product (i.e HHR norm).
	
	\begin{proposition}\label{normmod}
		Let $P \subset S(V)$ be a discrete $G$-subset. Then $P$ determines a natural module structure of $X$ over $N_*^PX$ as well as a unital $N_*^PX$-module map $N_*^PX \to X$.
		\begin{proof}
			Choose disjoint $G$-equivariant disks around each point $p$ in $S^V$ in a $G$-invariant way, and let $U$ be the union of these disks. For each orbit $O \subset P$ with isotropy group $H$, the union of the preimages of the disks around $O$ in $S_U$ is isomorphic to $G\times_H \RR^{V_{|H}}$. It then follows from properties $(F1)$ and $(F2)$ of equivariant factorization homology that $\int_{S_U}X\cong N_*^PX$. The result then follows by applying \Cref{mod} to $U$.
		\end{proof}
	\end{proposition}
	
	The following facts about $N^{\lambda}X$ are important to the proof of the main theorem:
	
	\begin{lemma}\label{underlyingdependence}
		Let $\lambda:G \to \Or(2)$ be a faithful representation. If $X\to Y$ is a map of $\EE_{\lambda}$-algebras that is an equivalence on $H$ fixed points for $H$ any proper subgroup of $G$, then the map $N^{\lambda}X\to N^{\lambda}Y$ is an equivalence of $\EE_{1}$-algebras.
		\begin{proof}
			In the case that $G \to \Or(2)$ factors through $\SOr(2)$, we can break up $\RR^{\lambda}-0$ into two copies of $\RR^2\times G$ glued along the collared manifold $(\RR^1\cup\RR^1)\times G\times (-1,1)$. By $G$-$\otimes$-excision $(F4)$ and the fact that induced manifolds are sent to norms $(F2)$, $\int_{\RR^2\times G}X = N_*^GX$, which depends on the underlying spectrum of $X$, giving the result.
			
			In the other case, $G=D_{2n}$ is a dihedral subgroup of $\Or(2)$ for $n\geq1$. If $n=1$, then $N^{\lambda}X$ = $N_*^{C_2}X$, which again only depends on the underlying spectrum.
			
			If $n\geq2$, we similarly use $G$-$\otimes$ excision to see that $N^{\lambda}X$ is identified the tensor product $N^{G/H_0}X\otimes_{N^GX}N^{G/H_1}X$ where $H_0$ are two conjugate maximal subgroups of $D_{2n}$. This formula again only depends on $X^H$ for proper subgroups of $G$.
		\end{proof}
	\end{lemma}
	
	\begin{remark}
		More generally, the factorization homology $\int_M$ depends on the $H$ fixed points when $H$ is contained in an isotropy group appearing in $M$. In particular $N^{V}X$ only depends on $X^{H}$ for $H$ a proper subgroup of $G$.
	\end{remark}
	
	When $\lambda$ is a two dimensional representation of $G$, $N^{\lambda}$ can be identified with a twisted version of $\THR$. We only use the following special case.
	\begin{lemma}\label{thhviathenorm}
		Let $\lambda:G \to \SOr(2)$ be a faithful representation, and let $X$ be a $G$-commutative algebra such that the $G$-action on the underlying algebra is trivial. Then there is an identification $N^{\lambda}X \simeq \THH(X_{un})$, where $\THH$ has the standard genuine $G$ action.
		\begin{proof}
			Consider the covering map $\RR^{\lambda}-0 \to \RR^2-0$ which (ignoring framings) is a $|G|$ fold cover. We can break $\RR^2-0$ into two copies of $\RR^2$ glued along a collared $(\RR\cup\RR)\times (-1,1)$, and this induces a decomposition of $\RR^{\lambda}-0$ into $\RR^2\times G$ and $(\RR\times \RR)\times G\times (-1,1)$. We can choose a framing on $(\RR\times \RR)\times G\times (-1,1)$ so that this gives a collar decomposition. The choice of framing doesn't matter because of commutativity.
			
			By $G$-$\otimes$-excision, this decomposition realizes $N^{\lambda}X$ as $N^GX\otimes_{N^GX\otimes N^GX}N^GX$. Since the action on the underlying spectrum of $X$ is trivial, this tensor product agrees with the formula in \cite[Theorem 4.4]{angeltveit2016topological} of the genuine $G$-action on $\THH(X)$.
		\end{proof}
	\end{lemma}
	
	In the next section we use:
	\begin{lemma}\label{pcomplete}
		Let $V$ be a representation of $G$ of dimension $\geq2$. Then $(N^{V}X_{V,p})^{\Phi G}$ and $X_{V,p}^{\Phi G}$ are $\FF_p$-modules.
		
		\begin{proof}
			Let $N$ be the kernel of $V$. Then $S(V)-S(V^G)$ contains an orbit of the form $G/N$. Thus we obtain by \Cref{normmod}, \Cref{snormmod}, and the naturality in \Cref{mod} an algebra map $N^{G}_NX_{V,p}\to N^V X_{V,p}$. We can apply geometric fixed points and use \cite[Proposition 2.57]{HHR} to obtain an algebra map $(X_{V,p})^{\Phi N} \to (N^{V}X_{V,p})^{\Phi G}$. But $N$ acts trivially on $V$ so it acts trivially on $\Omega^{V}S^{V+1}$ and we can identify $(X_{V,p})^{\Phi N}$ with the nonequivariant spectrum $X_{V,p}$. But this underlying spectrum is an $\FF_p$-algebra by \Cref{mahowald} and \Cref{hopkins} since $\dim V \geq 2$.
			
			Via \Cref{snormmod}, $X_{V,p}^{\Phi G}$ is a module over $(N^{V}X_{V,p})^{\Phi G}$, so it is also an $\FF_p$-module.
		\end{proof}
	\end{lemma}
	
	\section{Computing geometric fixed points}\label{geometricfixedpoints}
	Here, we study the geometric fixed points of $\ul{\FF}_p$ and $X_{\lambda,p}$, the objects of study in \Cref{mainthm}. Often $X_{\lambda,p}$ will be shortened to just $X_{\lambda}$. The results of this section together prove:
	
	\begin{customthm}{C}
		Let $\lambda:G \to \Or(2)$ be a faithful irreducible representation of a $p$-group. Then $X_{\lambda,p}^{\Phi G}$ and $\ul{\FF}_p^{\Phi G}$ are rank $2$ free modules over $(N^{\lambda}X_{\lambda,p})^{\Phi G}$ and $(N^{\lambda}\ul{\FF}_p)^{\Phi G}$ respectively, generated in degrees $0$ and $1$.
		\begin{proof}
			This follows immediately from \Cref{mackeyfreemodrk2}, \Cref{dihedralmackeyrank2}, and \Cref{freemodrk2}.
		\end{proof}
	\end{customthm}
	
	\subsection*{Thom spectra}
	In this subsection, we study $X_{\lambda}$ as a $N^{\lambda}X_{\lambda}$-module. The key fact allowing us to do so is that the module structure arises from taking the Thom spectrum of an action on the level of spaces.
	
	\begin{lemma}\label{1connected}
		Let $\lambda$ be a $2$-dimensional faithful irreducible representation of $G$.
		
		Then $\Map^G(S(\lambda),S^{\lambda+1})$ is simply connected.
		\begin{proof}
			For any such $\lambda$, $S(\lambda)$ has a cell decomposition with a $1$-cell with isotropy $G$ and zero cells with nontrivial isotropy group. Let $P$ be the $0$-skeleton of such a cell decomposition. Then there is a pushout square
			\begin{center}
				\begin{tikzcd}
					(S^0\times G)\ar[d]\ar[r] & P\ar[d]\\
					G\ar[r] & S(\lambda)\pushout
				\end{tikzcd}
			\end{center}
			which upon applying $\Map^G(-,S^{\lambda+1})$ yields the pullback square
			\begin{center}
				\begin{tikzcd}
					\Map^G(S(\lambda),S^{\lambda+1})\pullback\ar[d]\ar[r] & S^3\ar[d]\\
					\Map^G(P,S^{\lambda+1})\ar[r] & S^3\times S^3
				\end{tikzcd}
			\end{center}
			Since $P$ has no fixed points, $\Map^G(P,S^{\lambda+1})$ is a product of spheres of dimension $\geq 2$ so is simply connected. Applying the exact sequence on homotopy groups to the above pullback square, we obtain that $\Map^G(S(\lambda),S^{\lambda+1})$ is simply connected.
		\end{proof}
	\end{lemma}
	
	\begin{proposition}\label{freemodrk2}
		Let $\lambda$ be a rank 2 irreducible representation of $G$. Then $X_{\lambda,p}^{\Phi G}$ is a rank $2$ free module over $(N^{\lambda}X_{\lambda,p})^{\Phi G}$, generated in degrees $0$ and $1$.
		
		\begin{proof}
			For an $\EE_V$-space $X$ and a $V$-framed manifold $M$, equivariant nonabelian Poincaré duality says that $\int_MX = \Map_*(M^+,X)$, where $M^+$ is the $1$-point compactification of $M$ \cite[Theorem 2.2]{hahn2020equivariant}.
			
			This shows that $N^{\lambda}(\Omega^\lambda S^{\lambda+1})$ agrees with $\Map_*((\RR^{\lambda}-0)^+,S^{\lambda+1})$, where $(\RR^{\lambda}-0)^+$ can also be described as either the cofibre of the inclusion $S^0 \to S^\lambda$ or $\Sigma S(\lambda)_+$.
			
			Recall that by \Cref{snormmod}, $X_{\lambda}$ is a module over $N^{\lambda }X$. By the naturality of equivariant nonabelian Poincaré duality in the manifold $M$, the $N^{\lambda}(\Omega^\lambda S^{\lambda+1})$-module structure on $\Omega^{\lambda}S^{\lambda+1}$ comes from applying $\Map_*(-,S^{\lambda+1})$ on the natural coaction of the cogroup $(\RR^{\lambda}-0)^+=\Sigma S(\lambda)^+$ on $(\RR^{\lambda})^+$. This coaction arises from the coaction of the suspension on the cofibre in the cofibre sequence $$S(\lambda)^+ \xrightarrow{\ee}*^+\to (\RR^{\lambda})^+ \to (\RR^{\lambda}-0)^+$$ 
			Thus the action of $N^{\lambda}(\Omega^\lambda S^{\lambda+1})$ is the action of loop space on the fibre in the sequence
			
			$$N^{\lambda}(\Omega^{\lambda}S^{\lambda+1}) \to\Omega^{\lambda}S^{\lambda+1} \to \Map_*(S^0,S^{\lambda+1})\xrightarrow{\ee^*} \Map_*(S(\lambda)^+,S^{\lambda+1})  $$
			
			Let $R$ denote the $\EE_1$-ring $(N^{\lambda}X_{\lambda,p})^{\Phi G}$.
			The geometric fixed points of a Thom spectrum is given by the Thom spectrum of the fixed points, so the $R$-module structure of $X_{\lambda,p}$ on arises from applying the Thom spectrum to the $G$-fixed points of the action of the first two maps in the above fibre sequence..
			
			By \Cref{1connected}, $\ee^*$ is nulhomotopic after taking $G$-fixed points since $(S^{\lambda+1})^G = S^1$. It follows that $\Omega^{\lambda}S^{\lambda+1}$ is the product $(N^{\lambda}(\Omega^{\lambda}S^{\lambda+1}))^G\times S^{1}$, and $N^{\lambda}(\Omega^{\lambda}S^{\lambda+1})^G$ acts only on the first component.
			
			The Thom isomorphism in $\FF_p$-homology identifies the action on homology with the action of $N^\lambda X_{\lambda}$ on $X_{\lambda}$ in homology. But the action on $\FF_p$-homology is a rank $2$ free module generated in degrees $0$ and $1$ by the homology of $S^1$. By \Cref{pcomplete}, $X_{\lambda},N^{\lambda}X_{\lambda}$ are $\ul{\FF}_p$-modules, so it follows that $X_{\lambda}$ is a rank $2$ free module generated in degrees $0,1$.
		\end{proof}
	\end{proposition}
	
	\begin{remark}\label{remarkcareful}
		A more careful analysis in \Cref{1connected} shows that when $\lambda$ is a cyclic subgroup, $(\Omega^{\lambda}S^{\lambda+1})^G$ is equivalent to $\Omega S^3\times S^1\times \Omega^2S^3$. When $\lambda$ is a dihedral group $D_{2^n}$ for $n \geq 2$, it is equivalent to $\Omega S^2\times \Omega S^2\times S^1\times \Omega^2 S^3$. It is interesting that the proof of \Cref{freemodrk2} doesn't require explicitly identifying these spaces. We do carry this more careful analysis out for the representation $\rho$ in \Cref{rhoalgebra}. 
	\end{remark}
	
	We now identify $X_{\rho,2}^{\Phi C_2}$, where $\rho$ is the regular representation of $C_2$. Since an $\EE_{\rho}$-algebra is in particular an $\EE_1$-algebra as $1 \subset \rho$, and $(-)^{\Phi C_2}$ is symmetric monoidal, $X_{\rho,2}^{\Phi C_2}$ is an $\EE_1$-algebra.
	
	\begin{lemma}\label{rhoalgebra}
		$X_{\rho,2}^{\Phi C_2}$ is the $\EE_1$-algebra $\FF_2\otimes \Sigma^\infty_+\Omega S^2$.
		\begin{proof}
			$X_{\rho,2}^{\Phi {C_2}}$ is a Thom spectrum of $(\Omega^{\rho}S^{\rho+1})^{C_2}$. To compute these fixed points, we apply $\Map^{C_2}_*(-,S^{\rho+1})$ to the cofibre sequence $\Sigma(C_2)_+ \to S^1 \to S^{\rho} \to \Sigma^2(C_2)_+$ to get a fibre sequence of loop spaces
			$$\Omega^2 S^3 \to (\Omega^{\rho}S^{\rho+1})^{C_2} \to \Omega S^2 \to \Omega S^3$$
			Since the second map is a loop map, it is nulhomotopic, so $(\Omega^{\rho}S^{\rho+1})^{C_2} = \Omega (S^2 \times \Omega S^3)$ as a loop space, and $\Omega^2S^3$ acts trivially on the $\Omega S^2$ component. Because the Thom spectrum functor is $G$-symmetric monoidal and takes fixed points to geometric fixed points, the Thom spectrum of the composite $\Omega^2S^3 = \Map(\Sigma^2(C_2)_+,S^{\rho+1})^{C_2} \to  (\Omega^{\rho}S^{\rho+1})^{C_2} \to \BGL_1(\SP_{(2)})$ is $(N_*^{C_2}X_{\rho,2})^{\Phi {C_2}}$. By \cite[Proposition 2.57]{HHR}, this is the underlying spectrum of $X_{\rho,2}$, which is $X_{2,2} = \FF_2$ by \Cref{mahowald}, so it follows that $X_{\rho,2}^{\Phi {C_2}}$ is the $\EE_1$-$\FF_2$-algebra $\FF_2\otimes\Sigma^\infty_+ \Omega S^2$.
		\end{proof}
	\end{lemma}
	
	\subsection*{Mackey functors}\label{mackeysubsection}
	
	Now we turn to understanding $\ul{\FF}_p$. First we study $\ul{\FF}_p^{\Phi G}$ as a $N^{\lambda} \ul{\FF}_p^{\Phi G}$-module for the cyclic $p$-subgroups of $\SOr(2)$, using the identification of $N^{\lambda}\ul{\FF}_p$ with $\THH(\ul{\FF}_p)$.
	\begin{proposition}\label{mackeyfreemodrk2}
		Let $G$ be a $p$-subgroup of $\SOr(2)$. Then $\ul{\FF}_p^{\Phi G}$ is a rank $2$ free module over $(N^{\lambda} \ul{\FF}_p)^{\Phi G}$ generated in degrees $0$ and $1$.
		\begin{proof}
			Let $H$ be the maximal proper subgroup of $G$.
			We have a commutative square of ring maps
			\begin{center}
				\begin{tikzcd}
					(N^{\lambda}\ul{\FF}_p)^{\Phi G} \ar[r]\ar[d] & \ul{\FF}_p^{\Phi G}\ar[d]\\
					((N^{\lambda} \ul{\FF}_p)^{\Phi H})^{tG/H} \ar[r] & (\ul{\FF}_p^{\Phi H})^{tG/H}
				\end{tikzcd}
			\end{center}
			$N^{\lambda}\ul{\FF}_p$ can be identified with $\THH(\FF_p)$ as a $G$-spectrum by \Cref{thhviathenorm}. Since $\THH(\FF_p)$ is p-cyclotomic, the homotopy ring of $(N^{\lambda}\ul{\FF}_p)^{\Phi G} = \THH(\FF_p)^{\Phi G}$ agrees with that of $\THH(\FF_p)$, and so is a polynomial algebra $\FF_p[x]$ where $|x| = 2$. By \cite[Corollary IV.4.13]{nikolaus2018}, the element $x$ becomes a unit in $\pi_*\THH(\FF_p)^{tC_p}$.
			
			This map $\THH(\FF_p)^{\phi C_p} \to \THH(\FF_p)^{tC_p}$ is identified with the left vertical map in the diagram above, so we learn that $x \in \pi_2(N^{\lambda}\ul{\FF}_p^{\Phi G})$ is sent to a unit via that map. A unit cannot be in the kernel of any nonzero ring map, so the image of $x$ is not in the kernel of lower horizontal map. Thus $x$ has nonzero image in $(\ul{\FF}_p^{\Phi H})^{tG/H}$, and thus also in $\ul{\FF}_p^{\Phi G}$.

			The ring $\pi_*(\ul{\FF}_p^{\Phi G})$ is well known to be $\FF_2[t]$ for $|t|=1$ when $p=2$, and and $\FF_p[b]\otimes \Lambda(s)$ with $|b| = 2,|s|=1$ for $p>2$ (see for example \cite[Proposition 3.18]{HHR} or \cite[Lemma 8.2]{hahn2019eilenbergmaclane}). In either case, it is a rank $2$ free module over the subalgebra generated by any nonzero class in degree $2$, with module generators in degrees $0,1$. Thus the claim follows.
		\end{proof}
	\end{proposition}
	
	\begin{remark}
		It is possible to identify $N^{\lambda}\ul{\FF}_2$ for the dihedral groups $D_{2^n}$ with the genuine action of $D_{2^n}$ on $\THR(\ul{\FF}_2)$ (see \cite{dotto2017real, hesselholt2015real,quigley2019parametrized}). An approach similar to that in \Cref{mackeyfreemodrk2} should prove the analogous result for the dihedral groups.
	\end{remark}
	
	Next, we study $\ul{\FF}_p^{\Phi G}$ where $G$ is one of the dihedral groups $D_{2^n}$ for $n\geq2$. First we will review a few different ways of thinking about the isotropy separation sequence \cite[Section 2.5.2]{HHR}, which for a $G$-spectrum $X$ is the cofibre sequence
	
	\[X_{hP} \xrightarrow{\tr} X^G \to X^{\Phi G}\]
	where $X_{hP} = (X\otimes EP_+)^G$, with $EP$ the classifying space for proper subgroups of $G$, and $\tr$ the map induced from $EP_+ \to *_+=S^0$.
	
	Let $\Orb'_G$ be the category of nontrivial orbits of $G$. One model for $EP_+$ is the colimit $\colim_{G/H \in \Orb'_G}G/H_+$, showing that an alternative formula for $X_{hP}$ is $\colim_{G/H\in \Orb'_G}X^H$, and the natural map to $X^G$ is the colimit of the transfer maps. From this description, we see the functor $(-)^{\Phi G}$ factors through certain fixed point functors.
	
	\begin{lemma}\label{subgroupcofinal}
		Let $G$ be a group, and let $N$ be a normal subgroup contained in all the maximal proper subgroups of $G$. Then for any $G$-spectrum $X$, $X^{\Phi G} = (X^{N})^{\Phi(G/N)}$, where $X^N$ is given the structure of a genuine $G/N$ spectrum.
		\begin{proof}
			The category of orbits for proper subgroups of $G$ contains the category of orbits of proper subgroups of $G/N$ as the subcategory of orbits with isotropy containing $N$. We claim this inclusion is final, which implies the two isotropy separation sequences computing $X^{\Phi G}$ and $(X^{N})^{\Phi(G/N)}$ agree, giving the result.
			
			To see finality, we need to check for any orbit $G/H \in \Orb'_G$, the category maps to orbits $G/H'$ with $H' \supset N$ is nonempty and connnected. It is nonempty since $N$ contains all maximal proper subgroups. To see it is connected, any two maps $G/H \to G/H',G/H''$, choose basepoints on the $G$-sets so $H,H',H''$ are the isotropy groups of the basepoints. Then the projections $G/H \to G/H',G/H''$ factor through $G/N$ as pointed $G$-sets.
		\end{proof}
	\end{lemma}
	
	Note the above lemma is useful for understanding $\ul{\FF}_p^{\Phi G}$ since $\ul{\FF}_p^N$ agrees with $\ul{\FF}_p$ for the group $G/N$.
	
	Another model for $EP_+$ is $S(\infty V)_+$, where $V$ is any $G$-representation such that $V^G = 0$ and $V^H \neq 0$ for any proper subgroup of $G$. There is a cofibre sequence $S(\infty V)_+ \to S^0 \to S^{\infty V}$, showing that $X^{\Phi G}$ can also be described as $(X[a_V^{-1}])^G$ where $a_V$ is the element in $\pi_{-V}^G\SP$ arising from the inclusion $S^0 \to S^V$.
	
	As a warm up, we'll compute $\pi_*(\ul{\FF}_2^{\Phi G}))$ when $G$ is the group $C_{2^n}$. This is a well known computation, but we present it to demonstrate the strategy and set up notation for computing the geometric fixed points for the dihedral groups $D_{2^n}$.
	
	\begin{lemma}\label{cyclicmackeygeo}
		Let $G = C_{2^n}$. Then the homotopy ring $\pi_*(\ul{\FF}_2^{\Phi G}))$ is a polynomial algebra on a generator $u_{\sigma}$ in degree $1$.
		\begin{proof}
			Let $\sigma$ be the nontrivial $1$-dimensional representation of $G$, and let $b$ be a choice of cyclic generator of $G$. Since $\sigma^{H}$ is nontrivial for all proper subgroups of $G$, $\ul{\FF}_2^{\Phi G} = \ul{\FF}_2[a_\sigma^{-1}] = \colim_{n}\ul{\FF}_2\otimes S^{n\sigma}$. There is a minimal cell structure on $S^{n\sigma}$ with a cell called `$1$' in dimension $0$ with isotropy group $C_{2^n}$ and a cell called $x_i$ with isotropy group $C_{2^{n-1}}$ in dimension $i$ for $1 \leq i \leq n$. Taking the fixed points of the $\ul{\FF}_2$ cellular chain complex, we obtain a chain complex over $\FF_2$ generated by $1$ and $x_i+bx_i$, with zero differential, showing that $\pi_*\ul{\FF}_2\otimes S^{n\sigma}$ is $\FF_2$ in dimensions $0$ to $n$. Taking the limit as $n \to \infty$, this computes $\pi_*(\ul{\FF}_2^{\Phi G})$ additively.
			
			To obtain the multiplicative structure, we observe that there is another choice of cells on $S^{n\sigma}$ coming from the standard cell structure on $S^{\sigma}$ and the decomposition $S^{n\sigma} \cong (S^\sigma)^{\otimes n}$. The isomorphism can be chosen so that it sends $1 \mapsto 1^{\otimes n}$ and $x_i \mapsto 1^{\otimes i}\otimes x_1\otimes  (x_1+bx_1)^{\otimes{n-i-1}}$. From this, the multiplication sends $x_i \otimes x_j \to x_{i+j}$, so the description as a polynomial algebra follows.
		\end{proof}
	\end{lemma}
	
	Next, we study the dihedral groups $D_{2^n}$ for $n\geq 2$.
	
	\begin{proposition}\label{dihedralmackeygeo}
		Let $G = D_{2^n}$ for $n \geq 2$. Then there is an isomorphism of rings $$\pi_*(\ul{\FF}_2^{\Phi G}) =  \FF_2[u_{\sigma_1},u_{\sigma_2},u_{\sigma'}]/(u_{\sigma_1}u_{\sigma_2}+u_{\sigma_1}u_{\sigma'} + u_{\sigma_2}u_{\sigma'})$$
		\begin{proof}
			The intersection of the maximal proper subgroups of $D_{2^n}$ is a normal subgroup of index $4$ with quotient $D_4 = C_2\times C_2$. Thus by \Cref{subgroupcofinal}, we can assume $G = C_2\times C_2$. Suppose that $G$ is generated by $b_1,b_2$, let $\sigma_i$ be the $1$-dimensional representation of $G$ whose kernel is $b_i$, and let $\sigma'$ be $\sigma_1\otimes \sigma_2$.
			
			The element $\sigma_1+\sigma_2+\sigma'$ has no fixed points and its restriction to all proper subgroups is nontrivial. Thus $X^{\Phi G} = X[a_{\sigma_1+\sigma_2+\sigma'}^{-1}]$ for all $G$-spectra $X$. Moreover, $a_{\sigma_1+\sigma_2+\sigma'} = a_{\sigma_1}a_{\sigma_2}a_{\sigma'}$. 
			
			Thus $\ul{\FF}_2^{\Phi G} = \colim_{n,m,p \to \infty}\ul{\FF}_2\otimes S^{n\sigma_1+m\sigma_2+p\sigma'}$. We can take the minimal cell structures used for $S^{\infty \sigma_i}$ and $S^{\infty\sigma'}$ used in \Cref{cyclicmackeygeo} and tensor them together, where we use the names $x_i,y_i,z_i$ for the cells in $\sigma_1,\sigma_2,\sigma'$ respectively. In the $\ul{\FF}_2$ cellular chain complex, the elements $(x_1+b_2x_1)\otimes 1\otimes 1,1\otimes(y_1+b_1y_1)
			\otimes 1, 1\otimes 1\otimes (z_1+b_1z_1)$ represent the classes $u_{\sigma_1},u_{\sigma_2},u_{\sigma'}$ respectively. The products work just as in \Cref{cyclicmackeygeo} for the same reason.
			
			By applying the cellular differential to the chain $\sum_{g \in G}g(x_1\otimes y_1\otimes z_1)$, we obtain the relation $u_{\sigma_1}u_{\sigma_2}+u_{\sigma_1}u_{\sigma'} + u_{\sigma_2}u_{\sigma'}= 0$ in the cohomology of the cellular complex. It is straightforward to check that the classes $u_{\sigma_1},u_{\sigma_2},u_{\sigma'}$ generate the homology of the chain complex as an algebra, and that there are no more relations (for example by counting the dimension of the homology).
		\end{proof}
	\end{proposition}
	
	\begin{proposition}\label{dihedralmackeyrank2}
		Let $G = D_{2^n}$ for $n\geq2$, and $\lambda$ the standard inclusion to $\Or(2)$. $\ul{\FF}_2^{\Phi G}$ is a rank $2$ free module over $N^{\lambda}\ul{\FF}_2$ generated in degrees $0$ and $1$.
		\begin{proof}
			We first reduce to the case $n=2$. Present $G$ as $\langle a,b|a^2,abab,b^{2^{n-1}}\rangle$, and consider subgroups $H_0=\langle a\rangle,H_1 = \langle ab\rangle$.
			$G$-$\otimes$-excision (F4) lets us identify $N^{\lambda}\ul{\FF}_2$ with the tensor product $N^{G/H_0}\ul{\FF}_2\otimes_{N^G\ul{\FF}_2}N^{G/H_1}\ul{\FF}_2$, so applying $(-)^{\Phi G}$, we get $\ul{\FF}_2^{\Phi H_0}\otimes_{\FF_2} \ul{\FF}_2^{\Phi H_1}$. By factoring the functor $(-)^{\Phi G}$ as $((-)^{\Phi\langle b^2\rangle})^{\Phi G/\langle b^2 \rangle}$, we are reduced to the case $n =2$, so our group is $C_2\times C_2$.
			
			We compute the maps $\ul{\FF}_2^{\Phi{H_i}}\cong (N^{G/H_i}\ul{\FF}_2)^{\Phi G}  \xrightarrow{N^{\Phi G}} \ul{\FF}_2^{\Phi G}$ on $\pi_*$ for $i = 0,1$. By symmetry we can assume $i = 0$. The map $\ul{\FF}_2 \to \beta\ul{\FF}_2:=\Map(EG_+,\ul{\FF}_2)$ induces a commutative square:
			
			\begin{center}
				\begin{tikzcd}
					\ul{\FF}_2^{\Phi {H_0}} \ar[r,"N^{\Phi G}"]\ar[d] & \ul{\FF}_2^{\Phi  G} \ar[d]\\
					\ar[r,"N^{\tau G}"] \ul{\FF}_2^{\tau{H_0}} &  \ul{\FF}_2^{\tau G}
				\end{tikzcd}
			\end{center}
			
			The element $u_{\sigma}$ in $\pi_1(\ul{\FF}_2^{\Phi H_0})$ that gets sent to an invertible element in $\pi_1(\ul{\FF}_2^{\tau H_0})$. It follows that its image in $\ul{\FF}_2^{\tau G}$, a nonzero ring, cannot be $0$, so $u_{\sigma}$ has nonzero image in $\ul{\FF}_2^{\Phi  G}$ too. 
			
			To figure out what the image in $\ul{\FF}_2^{\Phi  G}$ is, we observe that $u_{\sigma}$ comes from applying geometric fixed points to an element in $\pi^{H_0}_{1-\sigma}\ul{\FF}_2$, where $\sigma$ is the nontrivial $1$-dimensional representation of $H_0$, so its norm is an element in $\pi^G_{\Ind_{H_0}^G(1-\sigma)}\ul{\FF}_2$. If $\sigma_1$ is the nontrivial representation of $G$ fixing $H_1$, $\sigma_0$ is the nontrivial representation fixing $H_0$, and $\sigma'$ is the tensor product, $\Ind_{H_0}^G(1-\sigma) = 1+\sigma_0 -\sigma_1-\sigma'$. Composing with $a_{\sigma_0}$, we get a nonzero element in $\pi^G_{1 -\sigma_1-\sigma'}\ul{\FF}_2$ in the image of the norm map. By the proof of \Cref{dihedralmackeygeo}, this element is a nonzero linear combination of $u_{\sigma'}$ and $u_{\sigma_1}$, but there is an automorphism fixing $H$ and swapping $\sigma'$ and $\sigma'$, so by naturality of $N^{\Phi G}$, it must be the sum $u_{\sigma_1}+u_{\sigma'}$. Symmetrically, the generator for $\pi_1(\ul{\FF}_2^{\Phi H_1})$ gets sent to $u_{\sigma_0}+u_{\sigma'}$.
			
			Finally, we observe that $\pi_*(\ul{\FF}_2^{\Phi G})) = \FF_2[u_{\sigma_0},u_{\sigma_1},u_{\sigma'}]/(u_{\sigma_1}u_{\sigma_0}+u_{\sigma_1}u_{\sigma'} + u_{\sigma_0}u_{\sigma'})$ is freely generated in degrees $0$ and $1$ as a module over $\FF_2[u_{\sigma_1}+u_{\sigma'},u_{\sigma_0}+u_{\sigma'}]$.
		\end{proof}
	\end{proposition}
	
	\section{Proof of Theorem A}
	
	In this section we prove \Cref{mainthm}. Let $\lambda$ be the standard $2$-dimensional representation of $\Or(2)$. We begin by producing a map of genuine $\Or(2)$-$\EE_{\lambda}$-algebras from $X_{\lambda,p} \to \ul{\FF}_p$. As before, $X_{\lambda,p}$ will often be shortened to $X_{\lambda}$.
	
	\begin{lemma}\label{mapproduce}
		There is a map $f_{\lambda}:X_{\lambda} \to \ul{\FF}_p$ of genuine $\Or(2)$-$\EE_{\lambda}$-algebras.
		\begin{proof}
			Since $X_{\lambda}$ is an $\EE_{\lambda}$ Thom spectrum, to give such a map $f_{\lambda}$, it suffices by \cite[X.6.4]{gaunce2006equivariant} to show the composite $S^{\lambda+1} \to \text{B}^{\lambda+1}\GL_1(\SP_{(p)}^0) \to \text{B}^{\lambda+1}\GL_1(\ul{\FF}_p)$ is null. 
			
			The map $\mu_{\lambda,p}$ represents the class $1-p$ in $\pi^G_0(\GL_1(\ul{\FF}_p))$, but $1-p = 1$ is the basepoint.
		\end{proof}
	\end{lemma}
	
	\begin{remark}
		In fact, the space of nulhomotopies in \Cref{mapproduce} is contractible, so the map $f_{\lambda}$ is essentially unique. 
	\end{remark}
	
	Restricting $f_{\lambda}$ to the $p$-subgroups of $\Or(2)$, we obtain our map of interest. There are two kinds of $p$-subgroups of $\Or(2)$. One family consists of the cyclic subgroups $C_{p^n} \subset \SOr(2)$, and the other family consists of the dihedral subgroups $D_{2^n},n\geq1$. By convention, $D_{2^n}$ has order $2^n$, so that $D_2$ is the cyclic group of order $2$ with its regular representation into $\Or(2)$.

	We prove \Cref{mainthm} by inducting on the finite subgroups of $\Or(2)$. This means the base case is when $G$ is the trivial group, where the result is the non-equivariant result of Hopkins--Mahowald (Theorems \ref{mahowald},\ref{hopkins}). This shows that Borel equivariantly, $f_{\lambda}$ is an equivalence, and the rest of the proof amounts to turning this into an equivalence of genuine $G$-spectra. In the inductive step, we have a finite $p$-subgroup $G \subset \Or(2)$, and we assume that $f_{\lambda}^{H}$ is an equivalence for all proper subgroups of $G$.
	
	The first step is to compute $\ul{\pi}_0$ of $X_{\lambda}$. The main ingredient is to use maps from norms of $X_{\lambda}$ to itself to show that certain elements of the Burnside ring are zero.
	
	\begin{proposition}\label{pi0}
		When $\lambda$ is a faithful rank $2$ representation of a $p$-group,	$\ul{\pi}_0(X_{\lambda}) = \ul{\FF}_p$.
		\begin{proof}
			$p=0$ in $\pi^G_0X_{\lambda}$ by construction: the Thom spectrum of the composite $S^1 \to \Omega^{\lambda}S^{\lambda+1} \to \BGL_1(\SP_{(p)})$ is the Moore spectrum $\SP/p$, so there is a unital map $\SP/p \to X_{\lambda}$, exhibiting the nulhomotopy of $p$ in $X_{\lambda}$.
			
			By induction on the subgroups of $G$, we can assume that the two Mackey functors $\ul{\FF}_p,\ul{\pi}_0X_{\lambda}$ agree on proper subgroups of $G$. The unit map $\SP \to X_{\lambda}$ is surjective on $\pi^G_0$ because $X_{\lambda}$ has only one cell in dimension $0$ and is connective, since $\Omega^{\lambda}S^{\lambda+1}$ is connected. As remarked already, $p=0$ in $\pi^G_0(X_{\lambda}^G)$. $\pi^G_0(\SP)$ is the Burnside ring $A(G)$, so it suffices to check that the classes $[G/H]$ are zero for every subgroup $H \lneq G$.
			
			If $H$ is a nonmaximal proper subgroup, then $[G/H]$ is zero by the induction hypothesis: we can choose $H \lneq H'\lneq G$, so that the class $[G/H]$ factors through the class $[H'/H]$, which we know to be zero in $X_{\lambda}^{H'}$.
			
			To finish the proof, we'll break into cases.
			
			If $G$ is the group $C_{p^n}$, any choice of a free orbit on $S(\lambda)$ gives by \Cref{normmod} a unital map $N_*^GX \to X$. Since $p = 0$, the composite map
			
			$$\SP = N_*^G\SP \xrightarrow{N_*^Gp} N_*^GX \to X$$
			is also zero. Since the map $N_*^GX \to X$ is unital, the diagram
			\begin{center}
				\begin{tikzcd}
					N_*^G\SP\ar[r,"N_*^Gp"] \ar[d,equal]& N_*^GX \ar[r] &X\\
					N_*^G\SP\ar[r,"N_*^Gp"] & N_*^G\SP\ar[r,equal] \ar[u,"N_*^G 1"]& \SP\ar[u,"1"]
				\end{tikzcd}
			\end{center}
			commutes. By the distributive law for the indexed tensor product \cite[A.3.3, Lemma A.36]{HHR}, the lower horizontal map corresponds in the Burnside ring to the $G$-set $\Hom(G,p)$. In the Burnside ring, $[\Hom(G,p)] \equiv [G/C_{p^{n-1}}] \pmod {p,[C_{p^{n-i}}],i\geq2}$ \cite[Lemma 6.6]{hahn2019eilenbergmaclane}, so it follows that $[G/C_{p^{n-1}}] = 0$.
			
			Next we consider the case $G=D_{2^n}$ for $n\geq2$. We present $D_{2^n}$ as $\langle a,b|a^2,b^{2^{n-1}},abab\rangle$ so that in the representation $\lambda$, $b \in \SOr(2)$, and $a$ acts as a reflection. There are three maximal proper subgroups of $G$: $\langle a,b^2\rangle, \langle ab,b^2\rangle, \langle b\rangle$. 
			
			$\RR^\lambda$ has orbits of the form $P_0 = G/\langle a\rangle,P_1 = G/\langle ab\rangle$, and $G$, so by \Cref{normmod} we obtain three unital maps from $N_*^GX_{\lambda}, N_*^{P_0}X_{\lambda}, N_*^{P_1}X_{\lambda}$ to $X$. As in the case of cyclic groups, this implies that $[\Hom(G,2)], [\Hom(P_0,2)],  [\Hom(P_1,2)]$ are all $0$. To see what this implies about $\pi^G_0(X_{\lambda})$, it suffices to compute each of these in the Burnside ring modulo the ideal generated by $2$ and $[G/H]$ for non-maximal proper subgroups $H$. 
			
			For each of the $G$-sets $S = G,P_0,P_1$, there are two fixed points in $\Hom(S,2)$. Since maximal subgroups of $G$ contain $\langle b^2\rangle$, we just need to count surjections in $\Hom(S,2)$ that are fixed by left multiplication by $\langle b^2\rangle$. 
			
			For $S = P_0$, it follows that $[\Hom(P_0,2)] \equiv [G/\langle a,b^2\rangle]$ since the double coset $\langle b^2\rangle \backslash G/\langle a\rangle$ consists of two elements. For the same reason, $[\Hom(P_1,2)] \equiv [G/\langle ab,b^2\rangle]$. For $S = G$, there is a contribution for each maximal proper subgroup, giving $[\Hom(G,2)] \equiv [G/\langle ab,b^2\rangle]+[G/\langle a,b^2\rangle]+[G/\langle b \rangle]$. Thus all of the desired classes in the Burnside ring are zero.
		\end{proof}
	\end{proposition}
	
	Below we compare the isotropy separation sequences for the map $f_{\lambda}:X_{\lambda} \to \ul{\FF}_p$.
	\begin{center}
		\begin{tikzcd}
			(X_{\lambda})_{hP} \ar[r]\ar[d,"(f_{\lambda})_{hP}"] & X_{\lambda}^G\ar[d,"f_{\lambda}^G"]\ar[r] & X_{\lambda}^{\Phi G}\ar[d,"f_{\lambda}^{\Phi G}"]\\
			(\ul{\FF}_p)_{hP} \ar[r] & \ul{\FF}_p^G \ar[r] &\ul{\FF}_p^{\Phi G}
		\end{tikzcd}
	\end{center}
	
	\begin{lemma}\label{surjpi1}
		Suppose $f_{\lambda}^H$ is an equivalence for all proper subgroups $H$ of $G$. Then the map $f_{\lambda}^{\Phi G}:X^{\Phi G}_{\lambda} \to \ul{\FF}_p^{\Phi G}$ is $1$-connected.
		\begin{proof}
			Since $f_{\lambda}^H$ is an equivalence for all proper subgroups of $G$, the map $(f_{\lambda})_{hP}$ is an equivalence, since $(f_{\lambda})_{hP} = \colim_{G/H\in \Orb'_G}f_{\lambda}^H$.
			
			By \Cref{pi0}, $f_{\lambda}^{G}$ is an isomorphism on $\pi_{\leq0}$. In addition, $\pi^G_1(\ul{\FF}_p) = 0$, so $f_{\lambda}^G$ is $1$-connected. Thus by the isotropy separation sequence, $f_{\lambda}^{\Phi G}$ is too.
		\end{proof}
	\end{lemma}
	
	We now finish the proof of \Cref{mainthm}.
	\begin{customthm}{A}
		Let $\lambda:G \to \Or(2)$ be a faithful representation of a $p$-group. Then there is an equivalence of $\EE_{\lambda}$-algebras $X_{\lambda,p}\cong \ul{\FF}_p$.
		\begin{proof}
			By induction, we can assume that the map $f_{\lambda}:X_{\lambda} \to \ul{\FF}_p$ induces an equivalence on fixed points for all proper subgroups. It suffices then to show $f_{\lambda}^{\Phi G}: X_{\lambda}^{\Phi G} \to \ul{\FF}_p^{\Phi G}$ is an equivalence. The base case is \Cref{mahowald} and \Cref{hopkins}.
			
			First consider the case when $\lambda$ is the regular representation $\rho$ of the group $C_2$. Then $X_{\lambda}^{\Phi C_2} \to \ul{\FF}_2^{\Phi C_2}$ is an $\EE_1$-algebra map. The codomain has homotopy ring $\FF_2[x], |x| = 1$, and the domain does too by \Cref{rhoalgebra}. By \Cref{surjpi1}, we learn that the generators are sent to the generators, so $f_{\lambda}^{\Phi C_2}$ is an equivalence.
			
			In all other cases, $\lambda$ is irreducible. For such $\lambda$,
			the induction hypothesis and \Cref{underlyingdependence} show that $N^{\lambda}X_{\lambda}^{\Phi G} \simeq N^{\lambda}\ul{\FF}_p^{\Phi G}$.

			Thus $f_{\lambda}^{\Phi G}$ is a map of $N^{\lambda}X_{\lambda}^{\Phi G} \simeq N^{\lambda}\ul{\FF}_p^{\Phi G}$-modules, and \Cref{freemodrank2both} shows that both sides are rank $2$ free modules generated in degrees $0$ and $1$. By \Cref{surjpi1}, both module generators of $\ul{\FF}_p^{\Phi G}$ are hit by $f_{\lambda}^{\Phi G}$. Thus the map is surjective on homotopy groups, but must also be injective since both sides are rank $2$ free modules.
		\end{proof}
	\end{customthm}

	\begin{remark}\label{evalgebradescription}
		The category $\Sp_G$ of $G$-spectra contains a full reflective subcategory $\Sp^P_G$ of $G$-spectra of the form $R^{P\beta}:=\ul{\Map}(EP,R)$, where $EP$ is the classifying space for proper subgroups of $G$. Any $\EE_V$-algebra in $\Sp_G$ $R$ can be presented by giving the data of $R^{P\beta}$ as an $\EE_V$-algebra in $\Sp^P_G$ along with an $\EE_{V^G}$-$N^VR$-algebra map $R^{\phi G} \to (R^{P\beta})^{\phi G}$.
		
		The proof of \Cref{mainthm} uses all this data, as we now explain. Fix $\lambda, G$ as in \Cref{mainthm}. In the proof, the inductive hypothesis gave that $(X_{\lambda,p})^{P\beta}\cong (\ul{\FF}_p)^{P\beta}$, and we used our understanding of $X_{\lambda,p}^{\phi G},\ul{\FF}_p^{\phi G}$ as $\EE_{\lambda^G}$-$N^{\lambda}\ul{\FF}_p^{\phi G}$-algebras to see that $X_{\lambda,p} \cong \ul{\FF}_p$.
	\end{remark}
	
	\begin{remark}
		The proof of \Cref{mainthm} in the case of the regular representation $\rho$ of $C_2$ is essentially the same as that in \cite{hahn2019eilenbergmaclane}. This is explained by the fact that $N^\rho = N_*^{C_2}$ is just the HHR norm. In all other cases overlapping with \Cref{hw}, the proof of \Cref{mainthm} improves the previously known proof, especially at odd primes, where an exotic Thom ring structure on $X_{\lambda,p}$, as well as an explicit computation of $\pi_1^G(X_{\lambda,p})$ is used in \cite{hahn2019eilenbergmaclane}.
	\end{remark}
	\section{Other Thom Spectra}\label{other}
	
	In this section we study other $X_{\lambda,p}$, showing that the choices of $G,\lambda,p$ in \Cref{mainthm} give all $X_{\lambda,p}$ isomorphic to $\ul{\FF}_p$. Then we discuss the algebraicity of $X_{\lambda,p}$ for other $G,\lambda,p$.

	\begin{customthm}{B}
			Let $G$ be a finite group and $V:G\to \Or(n)$ a representation. $X_{V,p} \cong \ul{\FF}_p$ if and only if $n = 2$, $G$ is a p-group, and $V$ is faithful.
		\begin{proof}
			The `if' part of the result is \Cref{mainthm}.
			
			The underlying spectrum of $X_{V,p}$ is a Thom spectrum of $\Omega^nS^{n+1}$ which only has the right mod $p$ homology if $\dim V = 2$, so that is a necessary condition.
			
			Note that if $X_{V,p}\cong \ul{\FF}_p$ holds for $V,p$, it must hold for $V_{|H},p$ for all subgroups $H$. Suppose some element $c$ of order $p$ acts trivially on $V$. Restricting to $H = \langle c \rangle$, we have the trivial representation on the cyclic group of order $p$. Since the action is trivial, there is an identification of $X_{V{|H},p}^{\Phi H}$ with the underlying Thom spectrum, which is $\FF_p$, and $\pi_*(\FF_p) \neq \pi_*(\ul{\FF}_p^{\Phi H})$, as the latter is $\FF_p$ in each nonnegative degree (see the discussion in \Cref{mackeyfreemodrk2}).
			
			Next, suppose that some element of prime order $q\neq p$ acts trivially on $V$. Let $H$ by the subgroup generated by this element so $V_{|H}$ is trivial. Then $X_{V_{|H},p} \ncong \ul{\FF}_p$ because $(\ul{\FF}_p)^{\Phi H} = 0$, but $X_{V_{|H},p}^{\Phi H} = \FF_p$ by Theorems \ref{hopkins}, \ref{mahowald}.
			
			Finally, suppose that $V$ is injective, and its image contains an element of prime order $q\neq p$. Then restricting to the subgroup generated by that element, \Cref{examplecomputation} (proven below) shows that $X_{V,p}$ doesn't agree with $\ul{\FF}_p$ when restricted to that subgroup. Thus the conditions on $V$ in \Cref{mainthm} are necessary.
		\end{proof}
	\end{customthm}
	
	For a finite group $G$ and an abelian group $A$, let $\tilde{A}$ be the Mackey functor whose value on $G/G$ is $A$ and otherwise $0$. Given a representation $V$, recall $a_V$ is the element in $\pi^G_{-V}\SP$ given by the inclusion $S^0 \to S^V$.
	
	\begin{theorem}\label{examplecomputation}
		Let $G=C_p$ be a cyclic group of prime order $p$, $\lambda:G \to \Or(2)$ an embedding, and $p\neq q$. Then as a $G$-spectrum, $X_{\lambda,q}\cong \ul{\FF}_q\oplus \bigoplus_0^\infty \Sigma^i\tilde{\FF}_ q$.
		\begin{proof}
			Let $X=X_{\lambda,q}$.	As a $C_p$ spectrum, there is a fracture square
			\begin{center}
				\begin{tikzcd}
					X\pullback\ar[r] \ar[d] & a_\lambda^{-1}X\ar[d]\\
					X_{a_{\lambda}}^{\wedge}\ar[r]& a_{\lambda}^{-1}(X_{a_{\lambda}}^{\wedge})
				\end{tikzcd}
			\end{center}
			which upon applying $G$-fixed points becomes the Tate square. By \Cref{pcomplete} and either \Cref{remarkcareful} or \cite[Proposition 3.9]{hahn2019eilenbergmaclane}, $X^{\Phi G} = a_\lambda^{-1}X^{G}$ is an $\FF_q$-module whose homotopy groups are $\FF_q$ in each nonnegative degree. Thus $a_{\lambda}^{-1}X^G \cong  \bigoplus_0^\infty \Sigma^i\tilde{\FF}_ q$. Furthermore, since the underlying action on $X$ is trivial on homotopy groups, and $p \neq q$, $X^{\wedge}_{a_{\lambda}}$ has fixed points equal to the homotopy fixed points, so is $\ul{\FF}_q$.
			
			The $G$-spectrum $a_{\lambda}^{-1}(X_{a_{\lambda}}^{\wedge})$ is zero, since its underlying spectrum is $0$ and the homotopy groups of its fixed points are the Tate cohomology of $\FF_q$ with a trivial $C_p$ action. Thus the above pullback square gives the decomposition claimed.
		\end{proof}
	\end{theorem}
	
	\begin{remark}
		Using \Cref{freemodrk2} and \Cref{evalgebradescription}, it is easy to see that the proof of \Cref{examplecomputation} gives a complete description of the $\EE_\lambda$-algebra structure.
	\end{remark}
	
	The $X_{\lambda,p}$ in \Cref{examplecomputation} are sums of generalized Eilenberg Mac Lane spectra, but this is not always the case. For example, take $\lambda$ to be a trivial $2$-dimensional representation for a nontrivial group, and let $p$ be any prime. Then $X_{2,p}$ is \textit{not} a generalized Eilenberg Mac Lane spectrum. Rather, it is the algebra $\SP\otimes \FF_p$, so that modules over it are the same as $\FF_p$-module valued Mackey functors. In this sense, it can still be considered algebraic. 
	
	In \Cref{conjecture}, it was asked what sufficient conditions are for $X_{\lambda,p}$ to be algebraic, i.e a $\FF_p$-module valued Mackey functor. It was also suggested that the only condition for $X_{\lambda,p}$ to be algebraic is that $\dim \lambda \geq 2$. The arguments of this paper are show that this is not far from being true. \Cref{pcomplete} showed that $X_{\lambda,p}^{\Phi G}$ is always an $\FF_p$-module in this case. By filtering $X_{\lambda,p}$ by its various isotropy separation sequences, we can use this to learn that $X_{\lambda,p}$ has a finite filtration whose associated graded objects are algebraic. Moreover, a refinement of the arguments of \Cref{examplecomputation} and \Cref{nogothm} show that $X_{\lambda,p}$ is algebraic whenever $G$ is a finite abelian group and $\dim \lambda \geq 2$.
	\section{Additional Structure}\label{deloop}
	
	 At the prime $2$, \Cref{mainthm} can be sharpened to an equivalence of $\EE_{\lambda+\det \lambda}$-algebras. To show this, one only needs to show that the map $\mu_{\lambda,2}:\Omega^{\lambda}S^{\lambda+1} \to \BGL_1\SP_{(2)}$ can be refined to a $\lambda+\det \lambda$-loop map. The goal of this section is to explain how to construct this refinement. We closely follow the strategy in \cite[Section 4]{hahn2019eilenbergmaclane}, except we bypass an unnecessary argument there. From now on, $\lambda$ denotes the standard $2$-dimensional representation of $\Or(2)$, and $\sigma$ denotes the representation $\det\lambda$.
	
	The ability to upgrade $\mu_{\lambda,2}$ to an $\EE_{\lambda+\sigma}$-loop map stems from the fact that $S^{1+\lambda}$ can be $\sigma$-delooped to $\HH\PP^\infty$, and equivariant Bott periodicity gives a natural $\sigma$-delooping of $\mu_{\lambda,2}$. First, we explain the delooping $\Omega^{\sigma}\HH\PP^\infty = S^{1+\lambda}$.
	
	The conjugation action of $\Sp(1)$ on $\HH$ makes $\HH$ into an $\SOr(3) \cong \Sp(1)/\{\pm1\}$-equivariant algebra, which induces an action of $\SOr(3)$ on $\HH\PP^\infty$. Restricting this action to $\Or(2)$, $\HH$ is equivariantly isomorphic to the Clifford algebra $\Cl(\lambda)$ where $\lambda$ is the standard representation of $\Or(2)$ with a negative definite form. As a representation, $\HH=\Cl(\lambda)$ is isomorphic to $1+\lambda+\det\lambda=1+ \lambda + \sigma$.
	
	Via the inclusion $\Or(2) \to \SOr(3)$, we can view $\HH\PP^\infty$ as a $\Or(2)$-space. Observe that $(\HH\PP^\infty)^G = \CC\PP^\infty$ when $G$ is a cyclic subgroup of $\Or(2)$, and that $(\HH\PP^\infty)^G = \RR\PP^\infty$ when $G$ is a dihedral subgroup $D_{2n}$ for $n\geq2$.

	\begin{lemma}
		Restricting the action of $\SOr(3)$ on $\HH\PP^\infty$ to $\Or(2)$, there is an equivalence of $\Or(2)$-spaces $S^{\lambda+1} \cong \Omega^{\sigma}\HH\PP^\infty$.
		\begin{proof}
			As an $\Or(2)$-space, $\HH\PP^1\cong S^{\HH} \cong S^{\lambda+\sigma+1}$, and so the inclusion $\HH\PP^1\to \HH\PP^\infty$ is adjoint to a map $S^{\lambda+1}\to \Omega^\sigma \HH\PP^\infty$. We show that this adjoint map is an equivalence on $G$ fixed points for every subgroup $G \subset \Or(2)$. For the cyclic subgroups of $\Or(2)$, this is \cite[Proposition 3.4,3.6]{hahn2019eilenbergmaclane}, so it suffices to show that that $S^1 \cong (\Omega^\sigma \HH\PP^\infty)^{D_{2n}}$ for $n\geq2$.
			
			There is a cofibre sequence
			$$(D_{2n}/C_{n})_+ \to S^0 \to S^{\sigma}$$
			Applying $\Map(-,\HH\PP^\infty)$, we obtain a fibre sequence
			$$\Omega^{\sigma}\HH\PP^\infty \to \RR\PP^\infty\to\CC\PP^\infty$$
			The fibre of the map $\text{B}C_2 \to \text{B}S^1$ is $S^1$, and which is $(S^{\lambda+\sigma+1})^{D_{2n}} = S^1$.
		\end{proof}
	\end{lemma}
	
	\begin{customthm}{D}
		Let $\lambda$ be the standard $2$-dimensional representation of $\Or(2)$. The map $\mu_{\lambda,2}:\Omega^{\lambda}S^{1+ \lambda} \to \BGL_1(\SP_{(2)})$ admits a refinement to a $\lambda+\det\lambda$-loop map.
		\begin{proof}
			The construction of our candidate refinement is exactly that in \cite[Section 4]{hahn2019eilenbergmaclane} which we refer to for more details. The difference between their arguments and what is done here is that we prove directly that the constructed map refines $\mu_{\lambda,2}$, bypassing the argument in their Proposition 4.12. 
			
			Let $\BGL(\HH)$ be the space classifying $G-\HH$-bundles, that is $G$-equivariant real vector bundles $E \to X$ with a $G$-equivariant algebra map $\HH \to \End(E)$. 
			
			There is an identification $\Omega^{\infty}\Sigma^{\HH}\KO_G\simeq \ZZ\times \BGL(\HH)$ \cite[Theorem 4.5]{hahn2019eilenbergmaclane}, where the equivalence is given by the Bott map: the tensor product of a bundle with the virtual $G-\HH$ bundle $\gamma$ on the pair $(\HH,\HH-\{0\})$ defined by the complex of bundles $\HH \to \HH$ given at a point $v\in\HH$ by multiplication by $v$ (see \cite[Section 4]{atiyah1966k}). Our delooping is the composite
			
			$$\HH\PP^\infty \xrightarrow{\cO(-1)-1} \BGL(\HH) \to \Omega^{\infty}\Sigma^{\HH} \KO_G \xrightarrow{\eta} \Omega^{\infty}\Sigma^{\lambda+\sigma} \KO_G $$ 
			
			To see that this does deloop $\mu_{\lambda,2}$, it suffices to identify the virtual $G-\HH$-bundle $\cO(-1)-1$ restricted to $\HH\PP^1$ with the class coming from $\gamma$. For then the element in $\pi^{\Or(2)}_{\lambda+\sigma+1}\Omega^{\infty}\Sigma^{\lambda+\sigma} \KO_G = \pi^{\Or(2)}_1 \BO_G$ is $\eta$, which is the class defining $\mu_{\lambda,2}$.
			
			To see this identification, we first analyze $\gamma$. Replacing the pair $(\HH,\HH-0)$ with the pair $(\HH,\HH-D^4)$ the data of the complex $\HH \to \HH$ defines a bundle $\gamma'$ over $Y=\HH_0\cup_{\HH-D^4}\HH_1$. There is a retraction $r:Y \to \HH_1$, and the corresponding virtual bundle is the difference of this bundle with the bundle on $Y$ that is the pullback of the composite $Y\xrightarrow{r}\HH_1 \to Y$, thought of as a relative class of the pair $(Y,\HH_1) \simeq (\HH,\HH-D^{\HH})$. However, this pullback is trivial since it factors through a contractible space, so after adding $1$ to the virtual bundle corresponding to $\gamma$, we just get the bundle $\gamma'$. There is an equivalence $Y \simeq \HH\PP^1$ identifying the copies of $\HH$ with hemispheres, under which $\gamma'$ corresponds to the bundle with the tautological clutching function $S(\HH) \to \GL_1(\HH)$, which is $\cO(-1)$.
		\end{proof}
	\end{customthm}
	\nocite{}
	\printbibliography
\end{document}